\documentclass{amsart}
\usepackage{amsmath}
\usepackage{amsfonts}
\usepackage{graphicx}
\usepackage{rotating}

\newtheorem*{theorem A}{Theorem A}
\newtheorem*{theorem B}{N\"olker's Theorem}

\theoremstyle{remark}

\theoremstyle{remark}

\theoremstyle{definition}

\numberwithin{equation}{section}
\def\({\left ( }
\def\){\right )}
\def\<{\left < }
\def\>{\right >}

\input{tcilatex}

\begin{document}
\title{Singular minimal translation graphs in Euclidean spaces}
\author{Ayla Erdur$^{1}$, Mahmut Ergut$^{2}$, Muhittin Evren Aydin$^{3}$ }
\address{$^{1,2}$ Department of Mathematics, Faculty of Science and Art,
Namik Kemal University, Tekirdag 59100, Turkey}
\address{$^{3}$ Department of Mathematics, Faculty of Science, Firat
University, Elazig, 23200, Turkey}
\email{meaydin@firat.edu.tr, aerdur@nku.edu.tr, mergut@nku.edu.tr}
\thanks{}

\begin{abstract}
In this paper, we consider the problem of finding the hypersurface $M^{n}$
in the Euclidean $\left( n+1\right) -$ space $\mathbb{R}^{n+1}$ that
satisfies an equation of mean curvature type, called singular minimal
hypersurface equation. Such an equation physically characterizes the
hypersurfaces in the upper halfspace $\mathbb{R}_{+}^{n+1}\left( \mathbf{u}%
\right) $ with lowest gravity center, for a fixed unit vector $\mathbf{u}\in 
\mathbb{R}^{n+1}$. We first state that a singular minimal cylinder $M^{n}$
in $\mathbb{R}^{n+1}$ is either a hyperplane or a $\alpha -$catenary
cylinder. It is also shown that this result remains true when $M^{n}$ is a
translation hypersuface and $\mathbf{u}$ a horizantal vector. As a further
application, we prove that a singular minimal translation graph in $\mathbb{R%
}^{3}$ of the form $z=f(x)+g(y+cx),$ $c\in \mathbb{R-\{}0\},$ with respect
to a certain horizantal vector $\mathbf{u}$ is either a plane or a $\alpha -$%
catenary cylinder.
\end{abstract}

\maketitle

\section{Introduction}

Let the pair $\left( \mathbb{R}^{3},g\right) $ denote the Euclidean $3-$%
space and $\mathbf{u}$ a fixed unit vector in $\mathbb{R}^{3}.$ Given a
smooth immersion $\sigma $ of an oriented surface $M^{2}$ into the halfspace 
\begin{equation*}
\mathbb{R}_{+}^{3}\left( \mathbf{u}\right) :\left\{ q\in \mathbb{R}%
^{3}:g\left( q,\mathbf{u}\right) >0\right\} .
\end{equation*}%
Let $\xi $ and $H$ be the Gauss map and the mean curvature of $\sigma $,
respectively (see \cite{Chen2}). Then, for some real constant $\alpha ,$ the
potential $\alpha -$energy of $\sigma $ in the direction of $\mathbf{u}$ can
be introduced by (see \cite{Dierkes1}, \cite{Lopez1}-\cite{Lopez3})

\begin{equation}
E\left( \sigma \right) =\underset{M^{2}}{\int }g\left( q,\mathbf{u}\right)
^{\alpha }dM^{2},  \tag{1.1}
\end{equation}%
where $dM^{2}$ denotes the measure on $M^{2}$ with respect to the induced
metric tensor from the Euclidean metric $g$ in $\mathbb{R}^{3}$ and $%
q=\sigma \left( p\right) ,$ $p\in M^{2}.$

Denoting $\Sigma :M^{2}\times \left( -\theta ,\theta \right) \rightarrow 
\mathbb{R}_{+}^{3}\left( \mathbf{u}\right) $ a compactly supported variation
of $\sigma $ with variaton vector field $\zeta ,$ the first variaton of $E$
then becomes 
\begin{equation*}
E^{\prime }\left( 0\right) =-\underset{M^{2}}{\int }\left( 2Hg\left( \sigma ,%
\mathbf{u}\right) -\alpha g\left( \xi ,\mathbf{u}\right) \right) g\left( \xi
,\zeta \right) ^{\alpha -1}dM^{2}.
\end{equation*}%
Taking $\sigma $ as a critical point of $E$, it then follows%
\begin{equation}
2H=\alpha \frac{g\left( \xi ,\mathbf{u}\right) }{g\left( \sigma ,\mathbf{u}%
\right) },  \tag{1.2}
\end{equation}%
which we call \textit{singular minimal }(\textit{SM})\textit{\ surface
equation}, see \cite{Dierkes2}. In the cited paper Dierkes calls a surface%
\textit{\ }$M^{2}$ a \textit{SM\ surface }or\textit{\ }$\alpha -$\textit{%
minimal surface\ }if its mean curvature $H$ fulfills Eq. (1.2).

Eq. (1.2) is clearly of mean curvature type and an extension of the classic
minimal surface equation (\cite[p. 17]{Lop}), i.e. the situation $\alpha =0$%
. Notice also that Eq. (1.2) is the Euler equation (see \cite[p. 33]{brunt})
of the variational integral given in Eq. (1.1).

Let $\gamma =\gamma \left( s\right) $ be a curve\textit{\ }in\textit{\ }$%
\mathbb{R}^{2}$ and $\mathbf{u}\in \mathbb{R}^{2}$ a fixed unit vector. For
the curve $\gamma $ one dimensional case of Eq. (1.2) writes%
\begin{equation}
\kappa \left( s\right) =\alpha \frac{g\left( N\left( s\right) ,\mathbf{u}%
\right) }{g\left( \gamma \left( s\right) ,\mathbf{u}\right) },  \tag{1.3}
\end{equation}%
where $\kappa $ and $N$ are the \textit{curvature} and unit \textit{%
principal normal vector field} of $\gamma .$ Hereinafter the curve $\gamma $
that the curvature $\kappa $ satisfies Eq. (1.3) is referred to as $\alpha -$%
\textit{catenary }\cite{Dierkes2}. Up to a change of coordinates, we can
take $\mathbf{u}=\left( 0,1\right) $ and $\gamma $ as the graph of $%
y=f\left( s\right) .$ In that case Eq. (1.3) writes%
\begin{equation}
\frac{f^{\prime \prime }}{1+f^{\prime 2}}=\frac{\alpha }{f}.  \tag{1.4}
\end{equation}%
Let the $y-$axis in\textit{\ }$\mathbb{R}^{2}$ denote the direction of the
gravity. What to solve Eq. (1.4) in case $\alpha =1$ is physically to find
the curve $\gamma $ in the upper halfplane $y>0$ with the lowest gravity
center \cite{Lopez4}. The solution of Eq. (1.4) is then the catenary 
\begin{equation*}
f\left( s\right) =\frac{1}{\lambda }\cosh \left( \lambda s+\mu \right)
,\lambda ,\mu \in \mathbb{R},\lambda \neq 0.
\end{equation*}

Let us now take the surface $M^{2}$ as a \textit{generalized cylinder }in%
\textit{\ }$\mathbb{R}^{3}$ (see \cite[p. 439]{Gray}) i.e. $M^{2}=\gamma
\left( s\right) +t\mathbf{w,}$ $s,t\in \mathbb{R}$, where $\gamma \left(
s\right) $ is the so-called \textit{base curve}, $\mathbf{w}\in \mathbb{R}%
^{3}$ a fixed unit vector. We call it \textit{cylinder}, shortly. Lopez \cite%
{Lopez1} proved that a SM cylinder is a $\alpha -$\textit{catenary cylinder, 
}namely a cylinder that takes the base curve as a $\alpha -$catenary. More
generally, in \cite{Lopez,Lopez1}, one was proved that a SM translation
surface, i.e. a graph of the form $z=f\left( x\right) +g\left( y\right) $ $,$
is a $\alpha -$catenary cylinder, where $\left\{ x,y,z\right\} $ is the
orthogonal coordinate system in $\mathbb{R}^{3}.$

In this paper, we generalize to higher dimensions the mentioned results in
previous paragraph. For this, we concern the following equation 
\begin{equation}
nH=\alpha \frac{g\left( \xi ,\mathbf{u}\right) }{g\left( \sigma ,\mathbf{u}%
\right) },\text{ }n\geq 2,  \tag{1.5}
\end{equation}%
where $\mathbf{u}\in \mathbb{R}^{n+1}$ is a fixed unit vector, $\xi $ and $H$
the Gauss map and the mean curvature of the smooth immersion $\sigma $ of an
oriented hypersurface $M^{n}$ into the halfspace 
\begin{equation*}
\mathbb{R}_{+}^{n+1}\left( \mathbf{u}\right) :\left\{ q\in \mathbb{R}%
^{n+1}:g\left( q,\mathbf{u}\right) >0\right\} .
\end{equation*}%
We call Eq. (1.5) \textit{SM\ hypersurface equation}.

Let $\left\{ x_{1},...,x_{n+1}\right\} $ denote the orthogonal coordinate
system in $\mathbb{R}^{n+1}$ and the $x_{n+1}-$axis the direction of the
gravity. Then Eq. (1.5) characterizes the hypersurfaces in $\mathbb{R}%
_{+}^{n+1}\left( \mathbf{u}\right) $ with lowest gravity center.

To study Eq. (1.5) we first take a generalized cylinder $M^{n}$ in $\mathbb{R%
}^{n+1}.$ By a \textit{generalized cylinder }in $\mathbb{R}^{n+1}$\textit{,}
we mean a hypersurface given in the following form%
\begin{equation}
M^{n}=\left\{ \gamma \left( s\right) +t_{1}\mathbf{w}_{1}+...+t_{n-1}\mathbf{%
w}_{n-1}:t_{1},...,t_{n-1}\in \mathbb{R},s\in I\subseteq \mathbb{R}\right\} ,
\tag{1.6}
\end{equation}%
where $\mathbf{w}_{1},\mathbf{w}_{2},...,\mathbf{w}_{n-1}$ are orthonormal
vectors in $\mathbb{R}^{n+1}$ and $\gamma $ (so-called \textit{base curve})
a 2-planar curve lying in $\Gamma =\mathtt{Span}\left\{ \mathbf{w}_{1},...,%
\mathbf{w}_{n-1}\right\} ^{\perp }.$ As before, we call it \textit{cylinder}%
, shortly. We prove that besides hyperplanes only SM cylinders in $\mathbb{R}%
^{n+1}$ are the $\alpha -$catenary cylinders.

Afterwards we study SM translation hypersurfaces in $\mathbb{R}^{n+1}$,
i.e., the graphs of the form (see \cite{Dillen})%
\begin{equation*}
x_{n+1}=f_{1}\left( x_{1}\right) +...+f_{n}\left( x_{n}\right) ,
\end{equation*}%
where $f_{1},...,f_{n}$ are smooth functions of single variable. We obtain
that besides hyperplanes SM translation hypersurfaces in $\mathbb{R}^{n+1}$
are $\alpha -$catenary cylinders. More geometric details on this class of
hypersurfaces can be found in \cite{Geomans}, \cite{HL}-\cite{H}, \cite%
{Lima,Li,Lop1}, \cite{Moruz}-\cite{Seo}.

As a further application, we concern SM translation graphs in $\mathbb{R}%
^{3} $ of the form $z=f(x)+g(y+cx),$ $c\in \mathbb{R-\{}0\}$. The study of
the surfaces of this kind was initiated by Liu and Yu \cite{Liu}, obtaining
explicit equations of minimal ones. These translation graphs, so-called 
\textit{affine translation surfaces,} belong to the family of surfaces
invariant by a group of translations (see Section 5)\textit{.} For some
results and progress on those surfaces, see \cite{Jung,Liu1}, \cite{Yang}-%
\cite{Yoon}. We show that a SM such translation graph in $\mathbb{R}^{3}$ is
either a plane or a $\alpha -$catenary cylinder.

\section{Preliminaries}

This section is devoted to present a short brief for hypersurfaces in $%
\mathbb{R}^{n+1}.$ Further details can be found in \cite{Chen3}.

Let $\sigma :M^{n}\rightarrow \mathbb{R}^{n+1}$ be a smooth immersion of an
oriented hypersurface $M^{n}$. Then the \textit{Gauss map} $\xi
:M^{n}\longrightarrow \mathbb{S}^{n}$ maps $M^{n}$ to the unit hypersphere $%
\mathbb{S}^{n}$ of $\mathbb{R}^{n+1}.$ The differential $d\xi $ of the Gauss
map $\xi $ is called \textit{Weingarten map}. For some vectors $\mathbf{v}$
and $\mathbf{w}$ tangent to $M^{n}$ at the point\ $p\in M^{n},$ the \textit{%
shape operator} $A_{p}$ is given by%
\begin{equation*}
\tilde{g}\left( A_{p}\left( \mathbf{v}\right) ,\mathbf{w}\right) =\tilde{g}%
\left( d\xi \left( \mathbf{v}\right) ,\mathbf{w}\right) ,
\end{equation*}%
where $\tilde{g}$ is the induced metric tensor on $M^{n}$ from the Euclidean
metric $g$ on $\mathbb{R}^{n+1}.$ The \textit{second fundamental form} $h$
of $\sigma $ is given in terms of the shape operator $A$ by%
\begin{equation*}
g\left( \xi ,h\left( \mathbf{v},\mathbf{w}\right) \right) =\tilde{g}\left(
A_{p}\left( \mathbf{v}\right) ,\mathbf{w}\right) .
\end{equation*}%
The \textit{mean curvature} of $\sigma $ at $p$ is defined by%
\begin{equation*}
H\left( p\right) =\frac{1}{n}\mathtt{tr}A_{p},
\end{equation*}%
where $\mathtt{tr}$ denotes the trace of $A_{p}.$ A hypersurface is called 
\textit{minimal} if $H$ vanishes identically.

The following result is well-known (see \cite{Chen1}).

\smallskip

\noindent \textbf{Proposition 2.1.} For a graph of $\mathbb{R}^{n+1}$ of the
form $x_{n+1}=f\left( x_{1},...,x_{n}\right) ,$ we have\newline
1. the unit normal vector field is 
\begin{equation*}
\xi =\frac{-1}{\phi }\left( f_{x_{1}},...,f_{x_{n}},-1\right) ,
\end{equation*}%
where $\phi =\sqrt{1+\sum_{j=1}^{n}\left( f_{x_{j}}\right) ^{2}}$ and $%
f_{x_{j}}=\frac{\partial f}{\partial x_{j}},$\newline
2. the components of the induced metric tensor (or the first fundamental
form) are%
\begin{equation*}
g_{ij}=\delta _{ij}+f_{x_{i}}f_{x_{j}},
\end{equation*}%
where $\delta _{ij}$ is Kronocker's Delta,\newline
3. the components of the second fundamental form are%
\begin{equation*}
h_{ij}=\frac{f_{x_{i}x_{j}}}{\phi },
\end{equation*}%
\newline
\newline
4. the matrix $\left[ a_{ij}\right] $ of the shape operator is%
\begin{equation*}
a_{ij}=\sum_{l}h_{il}g^{lj}=\frac{f_{x_{i}x_{j}}}{\phi }-\sum_{l}\frac{%
f_{x_{i}x_{l}}f_{x_{l}}f_{x_{j}}}{\phi ^{3}},
\end{equation*}%
where $f_{x_{i}x_{j}}=\frac{\partial ^{2}f}{\partial x_{i}\partial x_{j}}$
and $\left[ g^{lj}\right] =\left[ g_{lj}\right] ^{-1},$\newline
5. the mean curvature $H$ is%
\begin{equation*}
H=\frac{1}{n}\sum_{j=1}^{n}\frac{\partial }{\partial x_{j}}\left( \frac{%
f_{x_{j}}}{\phi }\right) .
\end{equation*}

\section{Generalized Cylinders}

Let $\mathbf{w}_{1},...,\mathbf{w}_{n-1}$ be fixed orthonormal vectors in $%
\mathbb{R}^{n+1}$ and $\gamma =\gamma \left( s\right) ,$ $s\in I\subseteq 
\mathbb{R}$ a unit speed curve lying in the $2-$plane $\Gamma =\mathtt{Span}%
\left\{ \mathbf{w}_{1},...,\mathbf{w}_{n-1}\right\} ^{\perp }.$ Consider the
cylinder in $\mathbb{R}^{n+1}$ given by Eq. (1.6). Denoting $\times $ the
cross product in $\mathbb{R}^{n+1}$ the unit normal vector field $\xi $ of $%
M^{n}$ becomes 
\begin{equation*}
\xi \left( s\right) =\mathbf{w}_{1}\times ...\times \mathbf{w}_{n-1}\times
\gamma ^{\prime }\left( s\right) =\sum_{i=1}^{n+1}\det \left( \mathbf{e}_{i},%
\mathbf{w}_{1},...,\mathbf{w}_{n-1},\gamma ^{\prime }\left( s\right) \right) 
\mathbf{e}_{i},
\end{equation*}%
where $\left\{ \mathbf{e}_{1},...,\mathbf{e}_{n+1}\right\} $ is the standard
basis of $\mathbb{R}^{n+1}$ and $\gamma ^{\prime }=\frac{d\gamma }{ds}$.
Then, as the components of the fundamental forms of $M^{n},$ we get $%
g_{ij}=\delta _{ij},$ $1\leq i,j\leq n,$ and%
\begin{equation*}
h_{ij}=\left\{ 
\begin{array}{l}
\kappa ,\text{ }for\text{ }i=j=n \\ 
0,\text{ }otherwise%
\end{array}%
\right.
\end{equation*}%
where $\kappa =\kappa \left( s\right) ,$ $s\in I,$ denotes the curvature of $%
\gamma .$ The mean curvature $H$ of $M^{n}$ turns to $H\left( s\right) =%
\frac{\kappa \left( s\right) }{n}$ and hence Eq. (1.5) leads to%
\begin{equation*}
\kappa \left( s\right) =\alpha \frac{g\left( \mathbf{w}_{1}\times ...\times 
\mathbf{w}_{n-1}\times \gamma ^{\prime }\left( s\right) ,\mathbf{u}\right) }{%
\sum_{i=1}^{n-1}g\left( \mathbf{w}_{i},\mathbf{u}\right) t_{i}+g\left(
\gamma \left( s\right) ,\mathbf{u}\right) }
\end{equation*}%
or equivalently%
\begin{equation}
\kappa \left( s\right) \sum_{i=1}^{n-1}g\left( \mathbf{w}_{i},\mathbf{u}%
\right) t_{i}+\kappa \left( s\right) g\left( \gamma \left( s\right) ,\mathbf{%
u}\right) -\alpha g\left( \mathbf{w}_{1}\times ...\times \mathbf{w}%
_{n-1}\times \gamma ^{\prime }\left( s\right) ,\mathbf{u}\right) =0. 
\tag{3.1}
\end{equation}%
If we take partial derivative of Eq. (3.1) with respect to $t_{i},$ $%
i=1,...,n-1,$ we find 
\begin{equation}
\kappa \left( s\right) g\left( \mathbf{w}_{i},\mathbf{u}\right) =0  \tag{3.2}
\end{equation}%
and%
\begin{equation}
\kappa \left( s\right) g\left( \gamma \left( s\right) ,\mathbf{u}\right)
-\alpha g\left( \mathbf{w}_{1}\times ...\times \mathbf{w}_{n-1}\times \gamma
^{\prime }\left( s\right) ,\mathbf{u}\right) =0.  \tag{3.3}
\end{equation}%
Assume, in Eq. (3.2), $\kappa =0$ on $I,$ identically. Then $\gamma $ is a
straight-line and Eq. (3.3) follows that $M^{n}$ is a hyperplane parallel to 
$\mathbf{u}$ because $\left\{ \mathbf{w}_{1},...,\mathbf{w}_{n-1},\gamma
^{\prime }\left( s\right) ,\mathbf{u}\right\} $ are linearly dependent at
every $s\in I$. Otherwise, i.e. $\kappa \neq 0,$ then $g\left( \mathbf{w}%
_{i},\mathbf{u}\right) =0$ and we get that $\mathbf{u}$ is parallel to $%
\Gamma .$ Hence Eq. (3.3) can be rewritten as%
\begin{equation}
\kappa \left( s\right) =\alpha \frac{g\left( \mathbf{w}_{1}\times ...\times 
\mathbf{w}_{n-1}\times \gamma ^{\prime }\left( s\right) ,\mathbf{u}\right) }{%
g\left( \gamma \left( s\right) ,\mathbf{u}\right) }=\alpha \frac{g\left(
N\left( s\right) ,\mathbf{u}\right) }{g\left( \gamma \left( s\right) ,%
\mathbf{u}\right) },  \tag{3.4}
\end{equation}%
where $N\left( s\right) $ denotes the principal unit normal vector field to $%
\gamma $ at $s\in I.$ Eq. (3.4) implies that $\gamma $ is a $\alpha -$%
catenary lying in $\Gamma $ and $M^{n}$ a $\alpha -$catenary cylinder.

Summing up, we have proved the following.

\smallskip

\noindent \textbf{Theorem 3.1.} \textit{The only SM cylinders in }$\mathbb{R}%
^{n+1}$\textit{\ are either the hyperplanes parallel to} $\mathbf{u}$ 
\textit{or the }$\alpha -$\textit{catenary\ cylinders whose the rulings are
orthogonal to }$\mathbf{u}$.

\section{Translation Hypersurfaces}

Let $\left\{ x_{1},...,x_{n+1}\right\} $ be the orthogonal coordinate system
in $\mathbb{R}^{n+1}.$ A \textit{translation hypersurface} $M^{n}$ in $%
\mathbb{R}^{n+1}$ can be described as the sum of $n$ curves $\gamma
_{1},...,\gamma _{n},$ so-called\textit{\ translating curves.} Then $M^{n}$
parameterizes locally as%
\begin{equation*}
\sigma \left( x_{1},...,x_{n}\right) =\gamma _{1}\left( x_{1}\right)
+...+\gamma _{n}\left( x_{n}\right) .
\end{equation*}%
If the translating curves $\gamma _{1},...,\gamma _{n}$ lie in orthogonal
2-planes mutually then after a change of coordinates $M^{n}$ becomes the
graph of the form%
\begin{equation*}
x_{n+1}=f_{1}\left( x_{1}\right) +...+f_{n}\left( x_{n}\right) ,
\end{equation*}%
where $f_{1},...,f_{n}$ are smooth functions of single variable. We mean
this graph by a translation hypersurface throughout the section.

By Proposition 2.1 the Gauss map $\xi $ and the mean curvature $H$ of $M^{n}$
are%
\begin{equation}
\xi =\frac{\left( -f_{1}^{\prime },...,-f_{n}^{\prime },1\right) }{\left[
1+\sum_{i=1}^{n}\left( f_{i}^{\prime }\right) ^{2}\right] ^{1/2}}\text{ } 
\tag{4.1}
\end{equation}%
and 
\begin{equation}
H=\frac{\sum_{i=1}^{n}\left( 1+\sum_{j\neq i}^{n}\left( f_{j}^{\prime
}\right) ^{2}\right) f_{i}^{\prime \prime }}{n\left[ 1+\sum_{i=1}^{n}\left(
f_{i}^{\prime }\right) ^{2}\right] ^{3/2}},  \tag{4.2}
\end{equation}%
where $f_{i}^{\prime }=\frac{df_{i}}{dx_{i}}$ and $f_{i}^{\prime \prime }=%
\frac{d^{2}f_{i}}{dx_{i}^{2}},$ $i=1,...,n.$ We have the following result:

\smallskip

\noindent \textbf{Theorem 4.1.} \textit{Let }$M^{n}$\textit{\ be a SM
translation hypersurface in} $\mathbb{R}^{n+1}$ \textit{with respect to a
horizantal vector} $\mathbf{u}$. \textit{Then it is either a hyperplane
parallel to }$\mathbf{u}$ \textit{or a }$\alpha -$\textit{catenary cylinder\
whose the rulings are horizontal straight-lines orthogonal to }$\mathbf{u}$.

\smallskip

\noindent \textbf{Proof. }Without lose of generality we can take $\mathbf{u}$
as $\left( 1,0,...,0\right) $. Eqs. (1.5), (4.1) and (4.2) then follow%
\begin{equation}
\sum_{i=1}^{n}\left( 1+\sum_{i\neq j}^{n}\left( f_{j}^{\prime }\right)
^{2}\right) f_{i}^{\prime \prime }=\frac{-\alpha f_{1}^{\prime }}{x_{1}}%
\left( 1+\sum_{i=1}^{n}\left( f_{i}^{\prime }\right) ^{2}\right) .  \tag{4.3}
\end{equation}%
We have to take $\alpha f_{1}^{\prime }\neq 0$ in Eq (4.3) because $M^{n}$
is minimal otherwise, which we are not interested in. Let us assume that $%
f_{1}^{\prime \prime }=0,$ or equivalently $f_{1}^{\prime }=const.\neq 0$,
in Eq. (4.3). Then the partial derivative of Eq. (4.3) with respect to $%
x_{1} $ leads to $\alpha =0$. This is the case we already ignore and
hereinafter $f_{1}^{\prime \prime }\neq 0$ is assumed. Next taking parital
derivative of Eq. (4.3) with respect to $x_{k},$ $k\neq 1,$ gives 
\begin{equation}
2f_{k}^{\prime }f_{k}^{\prime \prime }\sum_{k\neq i}^{n}f_{i}^{\prime \prime
}+\left( 1+\sum_{k\neq i}^{n}\left( f_{i}^{\prime }\right) ^{2}\right)
f_{k}^{\prime \prime \prime }=\frac{-2\alpha f_{1}^{\prime }}{x_{1}}%
f_{k}^{\prime }f_{k}^{\prime \prime }.  \tag{4.4}
\end{equation}%
It can be seen that $f_{k}^{\prime \prime }=0,$ $k=2,3,...,n,$ is a solution
of Eq. (4.4). This means that $f_{k}\left( x_{k}\right) =\lambda
_{k}x_{k}+\mu _{k},$ $\lambda _{k},\mu _{k}\in \mathbb{R},$ and $M^{n}$ is a
cylinder that can be written as%
\begin{eqnarray*}
\sigma \left( x_{1},x_{2},...,x_{n}\right) &=&\left(
x_{1},0,...,0,f_{1}\left( x_{1}\right) +\sum_{i=2}^{n}\mu _{k}\right)
+x_{2}\left( 0,1,...,0,\lambda _{2}\right) + \\
&&+...+x_{n}\left( 0,0,...,1,\lambda _{n}\right)
\end{eqnarray*}%
Due to Theorem 3.1, we obtain the the statement in the hypothesis of the
theorem. In order to finish the proof of the theorem, it is needed to show
that Eq. (4.4) has no other solution than $f_{k}^{\prime \prime }=0.$
Assuming now $f_{k}^{\prime \prime }\neq 0$ and dividing (4.4) with $%
2f_{k}^{\prime }f_{k}^{\prime \prime },$ we have%
\begin{equation}
\sum_{k\neq i}^{n}f_{i}^{\prime \prime }+\left( 1+\sum_{k\neq i}^{n}\left(
f_{i}^{\prime }\right) ^{2}\right) \frac{f_{k}^{\prime \prime \prime }}{%
2f_{k}^{\prime }f_{k}^{\prime \prime }}=\frac{-\alpha f_{1}^{\prime }}{x_{1}}%
.  \tag{4.5}
\end{equation}%
Taking parital derivative of Eq. (4.5) with respect to $x_{k},$ $k\neq 1,\ $%
gives%
\begin{equation*}
\left( 1+\sum_{k\neq i}^{n}\left( f_{i}^{\prime }\right) ^{2}\right) \left( 
\frac{f_{k}^{\prime \prime \prime }}{2f_{k}^{\prime }f_{k}^{\prime \prime }}%
\right) ^{\prime }=0,
\end{equation*}%
or equivalently,%
\begin{equation}
f_{k}^{\prime \prime \prime }=2\nu _{k}f_{k}^{\prime }f_{k}^{\prime \prime },
\tag{4.6}
\end{equation}%
for some constant $\nu _{k}.$ We distinguish two cases:

\begin{itemize}
\item $\nu _{k}=0,$ i.e. $f_{k}^{\prime \prime }=\lambda _{k},$ for some
nonzero constant $\lambda _{k},$ $k=2,3,...,n.$ After substituting this into
Eq. (4.3), one can be rewritten as%
\begin{equation}
\left. 
\begin{array}{c}
G_{1}\left( x_{1}\right) +G_{2}\left( x_{1}\right) \left( f_{2}^{\prime
}\right) ^{2}+G_{3}\left( x_{1}\right) \left( f_{3}^{\prime }\right)
^{2}+...+G_{n}\left( x_{1}\right) \left( f_{n}^{\prime }\right) ^{2}=0,%
\end{array}%
\right.  \tag{4.7}
\end{equation}%
where%
\begin{equation}
\left. 
\begin{array}{l}
G_{1}\left( x_{1}\right) =f_{1}^{\prime \prime }+\frac{\alpha f_{1}^{\prime }%
}{x_{1}}+\frac{\alpha \left( f_{1}^{\prime }\right) ^{3}}{x_{1}}+\left[
\left( f_{1}^{\prime }\right) ^{2}+1\right] \sum_{i=2}^{n}\lambda _{i} \\ 
G_{2}\left( x_{2}\right) =f_{1}^{\prime \prime }+\frac{\alpha f_{1}^{\prime }%
}{x_{1}}+\sum_{i=3}^{n}\lambda _{i} \\ 
G_{3}\left( x_{3}\right) =f_{1}^{\prime \prime }+\frac{\alpha f_{1}^{\prime }%
}{x_{1}}+\sum_{3\neq i=2}^{n}\lambda _{i} \\ 
\vdots \\ 
G_{n}\left( x_{1}\right) =f_{1}^{\prime \prime }+\frac{\alpha f_{1}^{\prime }%
}{x_{1}}+\sum_{i=2}^{n-1}\lambda _{i}.%
\end{array}%
\right.  \tag{4.8}
\end{equation}%
Because $f_{k}^{\prime \prime }\neq 0,$ $k=2,3,...,n,$ taking partial
derivative of (4.7) with respect to $x_{k}$ gives that the functions $%
G_{1},...,G_{n}$ are all zero. If we substract second equality in Eq. (4.8)
from third one then we find $\lambda _{2}=\lambda _{3}.$ Analogously if we
substract third equality in Eq. (4.8) from fourth one then we find $\lambda
_{3}=\lambda _{4}.$ Hence by proceeding this for other equalities in Eq.
(4.8) we obtain $\lambda _{2}=\lambda _{3}=...=\lambda _{n}.$ Put $\tilde{%
\lambda}=\lambda _{k},$ $k=2,3,...,n.$ The following can be obtained by some
equality in Eq. (4.8) (other than the first one) 
\begin{equation}
f_{1}^{\prime \prime }+\frac{\alpha f_{1}^{\prime }}{x_{1}}=\tilde{\lambda}%
(2-n).  \tag{4.9}
\end{equation}%
Substituting (4.9) into the first equality in Eq. (4.8) leads to%
\begin{equation}
\tilde{\lambda}(n-1)\left( f_{1}^{\prime }\right) ^{2}+\frac{\alpha \left(
f_{1}^{\prime }\right) ^{3}}{x_{1}}+\tilde{\lambda}=0.  \tag{4.10}
\end{equation}%
By taking derivative of Eq. (4.10) with respect to $x_{1}$ and then dividing 
$x_{1}$ we derive%
\begin{equation}
f_{1}^{\prime \prime }\left[ 2\tilde{\lambda}(n-1)\frac{f_{1}^{\prime }}{%
x_{1}}+3\alpha \left( \frac{f_{1}^{\prime }}{x_{1}}\right) ^{2}\right]
-\left( \frac{f_{1}^{\prime }}{x_{1}}\right) ^{3}=0.  \tag{4.11}
\end{equation}%
From Eq. (4.9) we get $f_{1}^{\prime \prime }=\tilde{\lambda}(2-n)-\frac{%
\alpha f_{1}^{\prime }}{x_{1}}$ and considering this into Eq. (4.11) yields
a polynomial equation of $\left( \frac{f_{1}^{\prime }}{x_{1}}\right) ,$ in
which the leading coefficient is $-3\alpha ^{2}-1.$ This leads to a
contradiction.

\item $\nu _{k}\neq 0.$ Hence Eq. (4.5) reduces to%
\begin{equation}
\sum_{k\neq i}^{n}f_{i}^{\prime \prime }+\nu _{k}\left( 1+\sum_{k\neq
i}^{n}\left( f_{i}^{\prime }\right) ^{2}\right) =\frac{-\alpha f_{1}^{\prime
}}{x_{1}}.  \tag{4.12}
\end{equation}%
Taking partial derivative of Eq. (4.12) with respect to $x_{l},$ $1\neq
l\neq k,$ yields%
\begin{equation}
f_{l}^{\prime \prime \prime }+2\nu _{k}f_{l}^{\prime }f_{l}^{\prime \prime
}=0.  \tag{4.13}
\end{equation}%
Because Eq. (4.6) hold for $k=2,3,...,n,$ we have $f_{l}^{\prime \prime
\prime }=2\nu _{l}f_{l}^{\prime }f_{l}^{\prime \prime }.$ By substituting
this into Eq. (4.13), we obtain 
\begin{equation}
\nu _{k}+\nu _{l}=0.  \tag{4.14}
\end{equation}%
On the other hand, integrating Eq. (4.6) gives%
\begin{equation}
f_{k}^{\prime \prime }=\nu _{k}\left( f_{k}^{\prime }\right) ^{2}+\mu _{k}. 
\tag{4.15}
\end{equation}%
Substituting Eqs. (4.14) and (4.15) into Eq. (4.12) leads to%
\begin{equation}
f_{1}^{\prime \prime }=-\left[ \nu _{k}\left( f_{1}^{\prime }\right) ^{2}+%
\frac{\alpha f_{1}^{\prime }}{x_{1}}+\varepsilon \right] ,  \tag{4.16}
\end{equation}%
where $\varepsilon =\nu _{k}+\sum_{k\neq i=2}^{n}\mu _{i}.$ After
substituting Eqs. (4.15) and (4.16) into Eq. (4.3) we can rearrange it as 
\begin{equation}
\left. 
\begin{array}{c}
\left( \sum_{i=2}^{n}\left[ \left( \nu _{i}-\nu _{k}\right) \left(
f_{i}^{\prime }\right) ^{2}+\mu _{i}\right] -\nu _{k}\right) \left(
f_{1}^{\prime }\right) ^{2}+ \\ 
+\sum_{i=2}^{n}\left( 1+\sum_{j\neq i=2}^{n}\left( f_{j}^{\prime }\right)
^{2}\right) \left( \nu _{i}\left( f_{i}^{\prime }\right) ^{2}+\mu
_{i}\right) + \\ 
-\varepsilon \left( 1+\sum_{i=2}^{n}\left( f_{i}^{\prime }\right)
^{2}\right) =\frac{-\alpha \left( f_{1}^{\prime }\right) ^{3}}{x_{1}}.%
\end{array}%
\right.  \tag{4.17}
\end{equation}%
The partial derivatives of Eq. (4.17) with respect to $x_{1}$ and $x_{l},$ $%
1\neq l\neq k$ lead to%
\begin{equation}
\nu _{k}-\nu _{l}=0.  \tag{4.18}
\end{equation}%
Comparing Eqs. (4.14) and (4.18) contradicts with $\nu _{k}\neq 0.$
\end{itemize}

\section{An application in 3-dimensional case}

Let $M^{2}$ be a translation graph in $\mathbb{R}^{3}$ of the form 
\begin{equation*}
z=f\left( x\right) +g\left( y+cx\right) ,\text{ }c\in \mathbb{R-\{}0\},
\end{equation*}%
for some smooth functions $f$ and $g.$ By the change of parameters $\tilde{x}%
=x$ and $\tilde{y}=y+cx$, we can choose a local parameterization on $M^{2}$
as 
\begin{equation*}
\sigma :I\times J\subset \mathbb{R}^{2}\rightarrow \mathbb{R}^{3}
\end{equation*}%
and%
\begin{equation}
\sigma \left( \tilde{x},\tilde{y}\right) =\left( \tilde{x},\tilde{y}-c\tilde{%
x},f\left( \tilde{x}\right) +g\left( \tilde{y}\right) \right) .  \tag{5.1}
\end{equation}%
Eq. (5.1) follows that it can be written as a sum of two planar curves, i.e. 
\begin{equation*}
M^{2}=\gamma \left( \tilde{x}\right) +\eta \left( \tilde{y}\right) ,
\end{equation*}%
where 
\begin{equation*}
\gamma :I\subset \mathbb{R}\rightarrow \mathbb{R}^{3},\text{ }\gamma \left( 
\tilde{x}\right) =\left( \tilde{x},-c\tilde{x},f\left( \tilde{x}\right)
\right) 
\end{equation*}%
and%
\begin{equation*}
\eta :J\subset \mathbb{R}\rightarrow \mathbb{R}^{3},\text{ }\eta \left( 
\tilde{y}\right) =\left( 0,\tilde{y},g\left( \tilde{y}\right) \right) .
\end{equation*}%
Notice that $\gamma $ and $\eta $  lie in the non-orthogonal planes to each
other. Hence $M^{2}$ turns to an extension of classical translation surface.

Let us put $f^{\prime }=\frac{df}{d\tilde{x}}$ and $g^{\prime }=\frac{dg}{d%
\tilde{y}},$ etc. Thereby we have\smallskip

\noindent \textbf{Theorem 5.1. }\textit{Let }$M^{2}$\textit{\ be a
translation graph in }$\mathbb{R}^{3}$ \textit{of the form }$z=f\left(
x\right) +g\left( y+cx\right) ,$\textit{\ }$c\neq 0.$\textit{\ If }$M^{2}$%
\textit{\ is a SM surface with respect to the horizantal vector} $\mathbf{u}%
=\left( 1,0,0\right) $\textit{\ then it is either a plane parallel to} $%
\mathbf{u}$ \textit{or} \textit{a }$\alpha -$\textit{catenary cylinder\
whose the rulings are horizontal straight-lines orthogonal to }$\mathbf{u}$%
.\smallskip

\noindent \textbf{Proof. }If $M^{2}$ is a SM surface then Eq. (1.2) writes%
\begin{equation}
\frac{\left[ 1+\left( g^{\prime }\right) ^{2}\right] f^{\prime \prime }+%
\left[ 1+c^{2}+\left( f^{\prime }\right) ^{2}\right] g^{\prime \prime }}{%
1+\left( f^{\prime }+cg^{\prime }\right) ^{2}+\left( g^{\prime }\right) ^{2}}%
=-\alpha \frac{f^{\prime }+cg^{\prime }}{\tilde{x}}.  \tag{5.2}
\end{equation}%
We distinguish several cases:

\begin{itemize}
\item $f^{\prime \prime }=0.$ Hence Eq. (5.2) reduces to 
\begin{equation}
\frac{\left[ 1+c^{2}+f_{0}^{2}\right] g^{\prime \prime }}{1+\left(
f_{0}+cg^{\prime }\right) ^{2}+\left( g^{\prime }\right) ^{2}}=-\alpha \frac{%
f_{0}+cg^{\prime }}{\tilde{x}},  \tag{5.3}
\end{equation}%
where $f^{\prime }=f_{0}$ is some constant. Because $\alpha \neq 0,$ the
partial derivative of Eq. (5.3) with respect to $\tilde{x}$ gives 
\begin{equation}
f_{0}+cg^{\prime }=0.  \tag{5.4}
\end{equation}
Eq. (5.4) implies $g^{\prime \prime }=0.$ Putting $g^{\prime }=g_{0},$ for
some nonzero constant $g_{0},$ we conclude%
\begin{equation*}
f_{0}+cg_{0}=0,
\end{equation*}%
which yields that 
\begin{equation*}
z\left( x,y\right) =f_{0}x+g_{0}\left( y+cx\right) +\mu =g_{0}y+\mu ,
\end{equation*}%
for some constant $\mu .$ This leads $M^{2}$ to a plane parallel to $\mathbf{%
u}$.

\item $f^{\prime \prime }\neq 0$ and $g^{\prime \prime }=0.$ Then we have $%
g\left( y+cx\right) =\lambda \left( cx+y\right) +\mu ,$ for some constants $%
\lambda ,\mu .$ Therefore $M^{2}$ turns to a cylinder of the form%
\begin{equation*}
\left. 
\begin{array}{c}
\sigma \left( x,y\right) =\left( x,y,f\left( x\right) +\lambda \left(
cx+y\right) +\mu \right) \\ 
=\left( x,0,f\left( x\right) +\lambda cx+\mu \right) +y\left( 0,1,\lambda
\right) .%
\end{array}%
\right.
\end{equation*}%
Due to \cite[Theorem 1]{Lopez1}, we obtain that $M^{2}$ is a $\alpha -$%
catenary cylinder whose the rulings are parallel to the vector $\left(
0,1,\lambda \right) $.

\item $f^{\prime \prime }g^{\prime \prime }\neq 0.$ By taking partial
derivative of Eq. (5.2) with respect to $\tilde{y},$ we have%
\begin{equation}
\left. 
\begin{array}{c}
2g^{\prime }g^{\prime \prime }f^{\prime \prime }+\left[ 1+c^{2}+\left(
f^{\prime }\right) ^{2}\right] g^{\prime \prime \prime }=\frac{-\alpha c}{%
\tilde{x}}\left[ 1+\left( f^{\prime }+cg^{\prime }\right) ^{2}+\left(
g^{\prime }\right) ^{2}\right] g^{\prime \prime }- \\ 
-2\alpha \frac{f^{\prime }+cg^{\prime }}{\tilde{x}}\left[ cf^{\prime
}+\left( c^{2}+1\right) g^{\prime }\right] g^{\prime \prime }.%
\end{array}%
\right.  \tag{5.5}
\end{equation}%
Dividing Eq. (5.5) with $g^{\prime \prime }$ leads to%
\begin{equation}
\left. 
\begin{array}{c}
2g^{\prime }f^{\prime \prime }+\left[ 1+c^{2}+\left( f^{\prime }\right) ^{2}%
\right] \frac{g^{\prime \prime \prime }}{g^{\prime \prime }}=\frac{-\alpha c%
}{\tilde{x}}\left[ 1+\left( f^{\prime }+cg^{\prime }\right) ^{2}+\left(
g^{\prime }\right) ^{2}\right] - \\ 
-2\alpha \frac{f^{\prime }+cg^{\prime }}{\tilde{x}}\left[ cf^{\prime
}+\left( c^{2}+1\right) g^{\prime }\right] .%
\end{array}%
\right.  \tag{5.6}
\end{equation}%
Assume now that $\frac{g^{\prime \prime \prime }}{g^{\prime \prime }}%
=\lambda ,$ for some constant $\lambda .$ Then Eq. (5.6) turns to%
\begin{equation}
\left. 
\begin{array}{c}
2\tilde{x}f^{\prime \prime }g^{\prime }+\left[ 1+c^{2}+\left( f^{\prime
}\right) ^{2}\right] \lambda \tilde{x}=-\alpha c\left[ 1+\left( f^{\prime
}+cg^{\prime }\right) ^{2}+\left( g^{\prime }\right) ^{2}\right] - \\ 
-2\alpha \left( f^{\prime }+cg^{\prime }\right) \left[ cf^{\prime }+\left(
c^{2}+1\right) g^{\prime }\right] .%
\end{array}%
\right.  \tag{5.7}
\end{equation}%
If we take partial derivative of Eq. (5.7) with respect to $\tilde{y}$ and
then divide it with $g^{\prime \prime },$ we find%
\begin{equation}
\tilde{x}f^{\prime \prime }+\alpha \left( 3c^{2}+1\right) f^{\prime
}=-3\alpha c\left( c^{2}+1\right) g^{\prime }.  \tag{5.8}
\end{equation}%
The partial derivative of Eq. (5.8) with respect to $\tilde{y}$ yields a
contradiction. Hence we conclude $\left( \frac{g^{\prime \prime \prime }}{%
g^{\prime \prime }}\right) ^{\prime }\neq 0.$ Next taking partial derivative
of Eq. (5.6) with respect to $\tilde{y}$ and dividing it with $g^{\prime
\prime }$ 
\begin{equation}
\left. 
\begin{array}{c}
f^{\prime \prime }+\frac{1}{2}\left[ 1+c^{2}+\left( f^{\prime }\right) ^{2}%
\right] \left[ \left( \frac{g^{\prime \prime \prime }}{g^{\prime \prime }}%
\right) ^{\prime }/g^{\prime \prime }\right] =\frac{-2\alpha c}{\tilde{x}}%
\left[ cf^{\prime }+\left( c^{2}+1\right) g^{\prime }\right] - \\ 
-\frac{\alpha \left( c^{2}+1\right) }{\tilde{x}}\left[ f^{\prime
}+cg^{\prime }\right] .%
\end{array}%
\right.  \tag{5.9}
\end{equation}%
The partial derivative of Eq. (5.9) with respect to $\tilde{y}$ yields%
\begin{equation}
\left[ 1+c^{2}+\left( f^{\prime }\right) ^{2}\right] \left[ \left( \frac{%
g^{\prime \prime \prime }}{g^{\prime \prime }}\right) ^{\prime }/g^{\prime
\prime }\right] ^{\prime }=\frac{-3\alpha c\left( c^{2}+1\right) }{\tilde{x}}%
g^{\prime \prime },  \tag{5.10}
\end{equation}%
which implies that both hand-sides cannot vanish. Hence Eq. (5.10) leads to 
\begin{equation}
1+c^{2}+\left( f^{\prime }\right) ^{2}=\frac{-3\alpha c\left( c^{2}+1\right) 
}{\lambda _{1}\tilde{x}}  \tag{5.11}
\end{equation}%
and 
\begin{equation}
\left[ \left( \frac{g^{\prime \prime \prime }}{g^{\prime \prime }}\right)
^{\prime }/g^{\prime \prime }\right] ^{\prime }=\lambda _{1}g^{\prime \prime
}  \tag{5.12}
\end{equation}%
for nonzero constant $\lambda _{1}.$ Integrating of Eq (5.12) gives%
\begin{equation}
g^{\prime \prime }=\frac{\lambda _{1}}{6}\left( g^{\prime }\right) ^{3}+%
\frac{\lambda _{2}}{2}\left( g^{\prime }\right) ^{2}+\lambda _{3}g^{\prime
}+\lambda _{4}  \tag{5.13}
\end{equation}%
for some constants $\lambda _{2},\lambda _{3},\lambda _{4}.$ Substituting
Eq. (5.13) into Eq. (5.2) gives%
\begin{equation}
\left. 
\begin{array}{c}
\left[ 1+\left( g^{\prime }\right) ^{2}\right] f^{\prime \prime }+\left[
1+c^{2}+\left( f^{\prime }\right) ^{2}\right] \left[ \frac{\lambda _{1}}{6}%
\left( g^{\prime }\right) ^{3}+\frac{\lambda _{2}}{2}\left( g^{\prime
}\right) ^{2}+\lambda _{3}g^{\prime }+\lambda _{4}\right] = \\ 
=-\alpha \frac{f^{\prime }+cg^{\prime }}{\tilde{x}}\left[ 1+\left( f^{\prime
}+cg^{\prime }\right) ^{2}+\left( g^{\prime }\right) ^{2}\right] .%
\end{array}%
\right.  \tag{5.14}
\end{equation}%
Eq. (5.14) is a polynomial equation of $g^{\prime }$ in which the
coefficient of the term of degree 1 satisfies%
\begin{equation}
\lambda _{3}\left[ 1+c^{2}+\left( f^{\prime }\right) ^{2}\right] +\frac{%
\alpha c}{\tilde{x}}+\frac{3\alpha c}{\tilde{x}}\left( f^{\prime }\right)
^{2}=0.  \tag{5.15}
\end{equation}%
Putting Eq. (5.11) into Eq. (5.15) yields%
\begin{equation*}
\frac{-3\alpha \lambda _{3}c\left( c^{2}+1\right) }{\lambda _{1}}+\alpha
c+3\alpha c\left( f^{\prime }\right) ^{2}=0,
\end{equation*}%
which implies $f^{\prime \prime }=0.$ This is not our case.
\end{itemize}

\section{Conclusions}

The results in Sections 4 and 5 on singular minimal translation
(hyper)surfaces were obtained by taking $\mathbf{u}$ as a horizontal vector
. It is pointed out that the vector $\mathbf{u}$ is parallel to the
hyperplane $x_{n+1}=0$ where the hypersurface is a graph. The investigating
of these graphs is still open when $\mathbf{u}$ is a vertical vector, that
is, $\mathbf{u}$ is normal to the hyperplane $x_{n+1}=0.$

Morever, the translation graph considering in Section 5 can be directly
extended to higher dimensions as%
\begin{equation}
z\left( x_{1},...,x_{n}\right) =f_{1}\left( x_{1}\right) +...f_{n-1}\left(
x_{n-1}\right) +g\left( x_{n}+\sum_{i=1}^{n-1}c_{i}x_{i}\right) ,\text{ }%
n\geq 2,  \tag{6.1}
\end{equation}%
for smooth functions $f_{1},...,f_{n-1},g.$ This is explicitly a
generalization in arbitrary dimensions of a classical translation
hypersurface. By taking the vector $\mathbf{u}$ as a horizontal or a
vertical vector, it could be another interesting problem to find a singular
minimal translation graph in $\mathbb{R}^{n+1}$ given by Eq. (6.1).

As a final remark these problems can be also considered by taking the vector 
$\mathbf{u}$ as arbitrary.

\end{document}